\documentclass[12pt]{article}
\usepackage{graphicx}
\usepackage[bf]{caption2}
\usepackage{amsfonts, amssymb, amsthm, amsmath}
\usepackage{color}
\newcounter{theorem}

\newcounter{theoremcounter}
\newcounter{lemmacounter}
\newcounter{remarkcounter}

\newtheorem{theorem}[theoremcounter]{Theorem}
\newtheorem{lemma}[lemmacounter]{Lemma}
\newtheorem{remark}[remarkcounter]{Remark}

\newcommand{\la}{\lambda}
\newcommand{\vp}{{\mathbf p}}
\newcommand{\vz}{{\mathbf z}}
\newcommand{\vf}{{\mathbf f}}
\newcommand{\vx}{{\mathbf x}}

\newcommand{\vy}{{\mathbf y}}
\newcommand{\vu}{{\mathbf u}}

\usepackage{enumerate}

\begin{document}

\title{Bounds on the Rate of Convergence for  $M^X_t/ M^X_t/1$ Queueing Models}

\author{Alexander Zeifman\thanks{Vologda State University, Vologda, Russia; Institute of Informatics Problems, Federal Research Center ``Computer Science and Control'',  Russian Academy of Sciences, Moscow, Russia; Vologda Research Center of the Russian Academy of Sciences, Vologda, Russia. E-mail: a$\_$zeifman@mail.ru} \and  Yacov Satin\thanks{Vologda State University, Vologda, Russia}\and  Alexander Sipin\thanks{Vologda State University, Vologda, Russia}}

\date{}

\maketitle

{\bf Abstract.} We apply the method of differential inequalities for 
 the computation of upper bounds for the rate
of convergence to the limiting regime for one specific class of
(in)homogeneous continuous-time Markov chains.  To obtain these estimates, we investigate the corresponding forward system of Kolmogorov differential equations.

\bigskip

{\bf Keywords:} {inhomogeneous continuous-time Markov chains, weak
ergodicity,  rate of convergence, sharp bounds, differential inequalities, forward Kolmogorov
system}

\section{Introduction}
\label{intro}

In this paper we consider the problem of finding the upper bounds
for the rate of convergence for some (in)homogeneous continuous-time Markov chains.

 To obtain these estimates, we investigate the corresponding forward system of Kolmogorov differential equations.

Consideration is given to classic inhomogeneous birth-death processes
and to special inhomogeneous chains with transitions intensities,
which do not depend on the current state.
Namely, let $\{X(t), \ t\geq 0 \}$ be an inhomogeneous continuous-time Markov chain
with the state space ${\mathcal{X}=\{ 0, 1, 2, \dots, \}}$.
Denote by $p_{ij}(s,t)=P\left\{ X(t)=j\left| X(s)=i\right. \right\}$,
$i,j \ge 0, \;0\leq s\leq t$, the transition probabilities of
$X(t)$ and by  $p_i(t)=P \left\{ X(t) =i \right\}$ -- the
probability that $X(t)$ is in state $i$ at time $t$.
Let $\vp(t) = \left(p_0(t), p_1(t), \dots,\right)^T$ be
probability distribution vector at instant $t$.
Throughout the paper it is assumed that in a small time interval
$h$ the possible transitions and their associated probabilities are
\begin{equation*}
p_{ij}(t,t+h)=
\left\{ \!
\begin{array}{cc}
q_{ij}(t)  h \!+\! \alpha_{ij}\left(t, h\right), & \mbox {if } j\neq i,
\\ 1 \! - \! \sum\limits_{k \in \mathcal{X}, k\neq i}q_{ik}(t) h \!+\! \alpha_{i}\left(
t,h\right), & \mbox {if } j=i,
\end{array}
\right.
%  \label{4001}
\end{equation*}
\noindent
where  $\sup_{i  \ge 0}  \sum_{j \ge 0}|\alpha_{ij}(t,h)| = o(h)$, for any $t \ge 0$. We also suppose that the transition intensities $q_{ij}(t) \ge 0$ are arbitrary non-random functions of $t$, locally integrable on $[0,\infty)$, and moreover, that there exists a positive number $L$ such that 
\begin{equation}
\sup_{i \in \mathcal{X}}
\left (\sum_{k \in \mathcal{X}, k\neq i}q_{ik}(t) \right )\le L < \infty,
\label{L}
\end{equation}
\noindent for almost all $t \ge 0$.
Then the probabilistic dynamics of the process $X(t)$ is
given by the forward Kolmogorov system
\begin{equation}
\label{ur01}
\frac{d}{dt}\vp(t)=A(t)\vp(t),
\end{equation}
\noindent where $A(t)$ is the transposed intensity matrix i.e.
$a_{ij}(t)=q_{ji}(t)$, $i,j \in \mathcal{X}$.  We can consider (\ref{ur01})  as the differential equation with bounded operator function in the space of sequences $l_1$ (see detailes, for instance in \cite{Zeifman2020dde}) and apply all results of \cite{DK}.

\smallskip

Throughout this paper by $\|\cdot\|$ (or by $\|\cdot\|_1$ if ambiguity is possible) we denote  the $l_1$-norm,
i.e. $\|{\vp(t)}\|=\sum_{i\in \mathcal{X}} |p_i(t)|$ and $\|A(t)\| =
\sup_{j \in \mathcal{X}} \sum_{i\in \mathcal{X}} |a_{ij}(t)|$. Let
$\Omega$ be a set of all stochastic vectors, i.e. $l_1$ vectors with
non-negative coordinates and unit norm. Then $\|A(t)\| \le 2L$ for almost all $t \ge 0$, and $\vp(s) \in \Omega$ implies 
$\vp(t) \in \Omega$ for any $0 \le s \le t$.

\smallskip

Recall that a Markov chain
$X(t)$ is called {\it weakly ergodic}, if
${\|\vp^{*}(t)-\vp^{**}(t)\| \to 0}$ as $t \to \infty$ for any
initial conditions $\vp^{*}(0)$ and $\vp^{**}(0)$, where
$\vp^{*}(t)$ and $\vp^{**}(t)$ are the corresponding solutions of
(\ref{ur01}).

\smallskip

We consider, as in \cite{Zeifman2020amcs}, the four classes of of Markov chains $X(t)$
with the following transition intensities:
\begin{enumerate}[i.]
  \item $q_{ij}(t)=0$ for any $t\ge 0$ if $|i-j|>1$ and both $q_{i,i+1}(t)=\la_i(t)$ and $q_{i,i-1}(t)=\mu_i(t)$ may depend on $i$;
  \item $q_{i,i-k}(t)=0$ for $k > 1$, $q_{i,i-1}(t)=\mu_i(t)$ may depend on $i$;  and $q_{i,i+k}(t)$, $k \ge 1$, depend only on $k$ and does not depend on $i$;
  \item $q_{i,i+k}(t)=0$ for $k > 1$, $q_{i,i+1}(t)=\la_i(t)$ may depend on $i$;  and $q_{i,i-k}(t)$, $k \ge 1$, depend only on $k$ and does not depend on $i$;
  \item both $q_{i,i-k}(t)$ and $q_{i,i+k}(t)$, $k \ge 1$, depend only on $k$ and do not depend on $i$.
\end{enumerate}

\smallskip

Each such process can be considered as the queue-length process for the corresponding queueing system $M^X_t/ M^X_t/1$.

Then type (i) transitions describe Markovian queues with possibly
state-dependent arrival and service intensities (for example, the
classic $M_n(t)/M_n(t)/1$ queue); type (ii) transitions allow
consideration of Markovian queues with state-independent batch
arrivals and state-dependent service intensity; type (iii)
transitions lead to Markovian queues with possible state-dependent
arrival intensity and state-independent batch service; type (iv)
transitions describe Markovian queues with state-independent batch
arrivals and batch service. We can refer to them as $M^X_t/ M^X_t/1$ queueing model following the original paper \cite{Nelson1988}, see also \cite{Li2017,Satin2013,Zeifman2014qs,Zeifman2020amcs}.

\smallskip

The simplest and most convenient for studying the rate of convergence to the limiting regime is the method of the logarithmic norm, see, for example \cite{Zeifman2014c,Zeifman2020dde,Zeifman2020amcs}. 

However, there are situations in which this approach does not give good results. 

Next, we show the possibility of using a different approach in such cases, namely, the method of differential inequalities.

Another (but similar) approach is to use piecewise-line Lyapunov functions, see, for example, \cite{Bertsimas2001,Blanchini2014,Bobyleva2002,Orlov2020}.

\smallskip

Due to the
normalization condition $p_0(t) = 1 - \sum_{i \ge 1} p_i(t)$, we
can rewrite the system (\ref{ur01}) as follows:
\begin{equation}
\frac{d}{dt}{\vz}(t)= B(t){\vz}(t)+{\vf}(t), \label{2.06}
\end{equation}
\noindent where
$$
{\vf}(t)=\left( a_{10}(t),  a_{20}(t),\dots \right)^{T}, \
{\vz}(t)=\left(p_1(t), p_2(t),\dots \right)^{T},
$$
\begin{equation}
{\footnotesize
\! B(t)\!=\!
%\!=\! \left(b_{ij}(t)\right)_{i,j=1}^{\infty} \!=\!
\left(
\begin{array}{ccccc}
a_{11} \!-\!a_{10} & a_{12}
\!-\!a_{10} & \cdots & a_{1r}
\!-\!a_{10} & \cdots \\
a_{21} \!-\!a_{20} & a_{22}
\!-\!a_{20} & \cdots & a_{2r}
\!-\!a_{20} & \cdots \\
\cdots & \cdots & \cdots & \cdots  & \cdots \\
a_{r1} \!-\!a_{r0} & a_{r2}
\!-\!a_{r0} & \cdots & a_{rr}
\!-\!a_{r0} & \cdots \\
\vdots & \vdots & \vdots & \vdots  & \ddots
\end{array}
\right)\!.}\label{2.07}
\end{equation}

\smallskip

Let  ${ \vy}(t) =
{\vz}^{*}(t)-{\vz}^{**}(t)$ be the difference of two solutions of
system (\ref{2.06}), and  ${ \vy}(t) = \left({y}_1(t), {y}_2(t),
\dots, \right)^T$.  Then, in contrast to the coordinates of the vector $\vp(t)$, 
the coordinates of the vector ${\vy}(t)$ have arbitrary signs.

\smallskip

Consider now  the  'homogeneous' system
\begin{equation}
\frac{d}{dt}{\vy}(t)= B(t){ \vy}(t), \label{hom1}
\end{equation}
\noindent  corresponding to (\ref{2.06}). As it was firstly noticed in \cite{Zeifman1989}, it
is more convenient to study the rate of convergence using the
transformed version $B^*(t)$ of $B(t)$ given by
$B^*(t)=TB(t)T^{-1}$, where $T$ is the upper triangular
matrix of the form
\begin{equation}
T=\left(
\begin{array}{cccccccc}
1   & 1 & 1 & \cdots & 1 & \cdots \\
0   & 1  & 1  &   \cdots & 1 & \cdots\\
0   & 0  & 1  &   \cdots & 1 & \cdots\\
\vdots & \vdots & \vdots & \ddots & \cdots\\
0   & 0  & 0  &   \cdots & 1 & \cdots \\
\cdots & \cdots & \cdots & \cdots & \cdots& \cdots
\end{array}
\right). \label{vspmatr}
\end{equation}

\noindent Let ${\vu}(t)=T{\vy}(t)$. Then the system \eqref{hom1} can be rewritten in the form
\begin{equation}
\frac{d}{dt} {\vu}(t) = B^*(t) {\vu}(t), \label{hom11}
\end{equation}

\noindent where ${\vu}(t) = \left( {u}_1(t), u_2(t), \dots\right)^T$
is the vector with the coordinates of arbitrary signs. If one of the two matrices  $B^*(t)$ or $B(t)$
is known, the other is also (uniquely) defined.

The approach based on the differential inequalities (see \cite{Zeifman2020amcs}) seems to be
the most general. On the other hand, if $ B^*(t) $ is essentially non-negative (i.e. all off-diagonal elements are non-negative for any $t \ge 0$), then the method based on the logarithmic norm gives the same results, but in a much more visual form, see \cite{Zeifman2020amcs}.

\section{Explicit forms of the reduced intensity matrices}

Here we only write out the form of the matrix $B^*(t)$ in each case, detailed obtaining of the corresponding transformations one can see in \cite{Zeifman2020amcs}.

\smallskip

For $X(t)$ belonging to class $\rm (i)$ (inhomogeneous birth-death process)
one has 
$$B^*(t) =  $$
\begin{equation}
{\tiny \left(
\begin{array}{ccccccc}
-\left(\lambda_0+\mu_1\right)  & \mu_1
 & 0 & \cdots & 0 &\cdots &\cdots\\
\lambda_1  & -\left(\lambda_1+\mu_2\right) & \mu_2 & \cdots & 0 &\cdots &\cdots\\
\ddots & \ddots & \ddots & \ddots & \ddots  &\cdots\\
0 & \cdots & \cdots & \lambda_{r-1} &
-\left(\lambda_{r-1}+\mu_r\right) & \mu_{r} &\cdots  \\
\cdots & \cdots & \cdots & \cdots & \cdots  & \cdots &\cdots \\
\end{array}
\right).} \label{class1-inf}
\end{equation}

\smallskip

For $X(t)$ belonging to class $\rm (ii)$ (which corresponds to the queueing system with batch arrivals and single services)
one has 
\begin{equation}
{\scriptsize B^*(t) =  \left(
\begin{array}{cccccc}
a_{11}  & \mu_1
 & 0 & \cdots & 0 \\
a_1  & a_{22} & \mu_2 & \cdots & 0 \\
a_2  & a_1  &a_{33} & \mu_3 & \cdots &  \\
\ddots & \ddots & \ddots & \ddots & \ddots  \\
\ddots & \ddots & \ddots & \ddots & \ddots  \\
\end{array}
\right)}. \label{class2-inf}
\end{equation}
\noindent

\smallskip

For $X(t)$ belonging to class $\rm (iii)$ (which corresponds to the queueing system with single arrivals and group services)
one has 
$B^*(t)=$
\begin{equation}
{\tiny \left(
\begin{array}{cccccc}
-\left(\lambda_0 +b_1\right)  & b_1 - b_2
 & b_2 - b_3 & \cdots & \cdots \\
\lambda_1 & -\big(\lambda_1+\sum\limits_{i\le 2}b_i\big) & b_1 - b_3 & \cdots & \cdots\\
\ddots & \ddots & \ddots & \ddots & \ddots  \\
0 & \cdots & \cdots & \lambda_{r-1} &
-\big(\lambda_{r-1}+\sum\limits_{i\le r}b_i\big) \cdots \\
\ddots & \ddots & \ddots & \ddots & \ddots  \\
\end{array}
\right)}.\label{class3-inf}
\end{equation}

\smallskip

Finally, for $X(t)$ belonging to class $\rm (iv)$ (which corresponds to the queueing system with state-independent batch arrivals and group services) one has

\begin{equation}
B^* = {\small \left(
\begin{array}{cccccc}
a_{11}  &  b_1 - b_2
 & b_2 - b_3 & \cdots & \cdots \\
a_1  & a_{22} & b_1 - b_3 & \cdots & \cdots \\ \\
\ddots & \ddots & \ddots & \ddots & \ddots  \\
a_{r-1} & \cdots & \cdots & a_1 & a_{rr} & \cdots \\
 \cdots & \cdots & \cdots & \cdots & \cdots & \cdots \\
\end{array}
\right)}.\label{class4-inf}
\end{equation}

\begin{remark} Generally speaking, for models of the first and second classes the matrix $B^*(t)$ is always essentially non-negative;
at the same time, for models of the third and fourth classes, this requires some additional assumptions.
Under essential non-negativity of $B^*(t)$ all bounds on the rate of convergence can be obtained via logarithmic norm, see \cite{Zeifman2020amcs}. However, in the general case, this approach may not work, and the method of differential inequalities described in our previous papers, see  \cite{Kryukova2020,Zeifman2020amcs} would be more effective.
\end{remark}

\smallskip

As a result, in this note we will consider chains of the third and fourth classes with countable state space. For simplicity of calculations, we will additionally assume that the size of the simultaneously arriving and/or servicing group of customers does not exceed some fixed number, say $R$, i.e. that all $q_{ij}(t) =0$ for $|i-j| >R$ and any $t \ge 0$.

\medskip

Let ${\{d_i, \ i \ge 1\}}$ be a sequence of non-zero
numbers such that $\inf_k |d_k| = d > 0$. Denote by  $D=diag(d_1,d_2,\dots)$  the corresponding diagonal matrix, with
the off-diagonal elements equal to zero. By putting ${\bf w}(t)=D
{\bf u}(t)$ in (\ref{hom11}), we obtain the following equation
\begin{equation}
\frac{d }{dt} {\bf w}(t)= B^{**}(t){\bf w}(t), \label{hom111}
\end{equation}
\noindent where 
\begin{equation}
B^{**}(t)=D
B^{*}(t)D^{-1}=\left(b^{**}_{ij}(t)\right)_{i,j \ge 1}.
\label{b**01}
\end{equation}

\smallskip

If we  write out  $B^{*}(t)= \left(b^{*}_{ij}(t)\right)_{i,j \ge1}$, then
\begin{equation}
b^{**}_{ij}(t)= \frac{d_i}{d_j}b^{*}_{ij}(t), \quad |i-j| \le R,
\label{b**02}
\end{equation}
\noindent and our assumption implies $b^{**}_{ij}(t) = b^{*}_{ij}(t) = 0$ for any $t \ge 0$ if $|i-j| > R$.

\section{Upper bounds on the rate of convergence}

Let us first consider a general {\it finite} system of linear differential equations, which we will write in the form
\begin{equation}
 \frac{d}{dt}\vx(t)=B^*(t)\vx(t), \ t \ge 0,
\label{upper31}
\end{equation}
\noindent where ${ \vx}(t) = \left({x}_1(t), \dots, {x}_S(t)\right)^T$, and let $D$ now be the corresponding {\it finite} diagonal matrix.

\smallskip

Firstly we consider a situation where all coefficients $b^*_{ij}(t)$ are analytical functions of $t$, it was considered in \cite{Kryukova2020,Zeifman2019icumt,Zeifman2020amcs}. The proof is based on the fact that in this case, on any finite interval, each coordinate has a finite number of sign changes, which means that the semiaxis can be divided into intervals, on each of which the signs of the coordinates are constant. Let, for instance at interval $(t_1,t_2)$ some coordinates of the solution are positive, and all others are negative.  Choose the signs of $d_k$-s  so that all $d_kx_k(t)>0$.
Hence $\|{\bf w}(t)\| = \|{\bf x}(t)\|_{D} = \sum_{k=1}^S d_kx_k(t) \ge d \|{\bf x}(t)\|_1$ can be considered as the corresponding norm.

Let $\sum_{i=1}^S b^{**}_{ij}(t) \le -\alpha_D(t)$, for any $j$,
then
\begin{equation}
\frac{d}{dt} \|{\bf w}(t)\|=
\frac{d\left(\sum_kw_k\right)}{dt}=
\sum_{i,j}b^{**}_{ij}(t)w_j(t)
\le -\alpha_D(t)\|{\bf w}(t)\|. \label{upper001}
\end{equation}
\noindent

\smallskip

Then we obtain the corresponding inequality
\begin{equation}
\|D{\bf x}\left(t\right)\|_1 \le e^{-\int_s^t\alpha_D (\tau) d\tau}\|D{\bf x}(s)\|_1, \quad t_1 < s < t <t_2,
\label{upper032}
\end{equation}
\noindent for the corresponding matrix $D$ and corresponding function $\alpha_D(t)$. 
\noindent Hence we have 
\begin{equation}
\|{\bf x}\left(t\right)\|_1 \le \frac{\max |d_k|}{\min |d_m|}e^{-\int_s^t\alpha_D (\tau) d\tau}\|{\bf x}(s)\|_1,
\label{upper033}
\end{equation}
\noindent for any $t_1 < s < t <t_2$. 

\smallskip

Then we glue such estimates on all intervals and get the general estimate 
\begin{equation}
\|{\bf x}\left(t\right)\|_1 \le d^*(S) e^{-\int_s^t\alpha^* (\tau) d\tau}\|{\bf x}(s)\|_1,
\label{upper034}
\end{equation}
\noindent where $\alpha^*(t) = \min \alpha_D(t)$, and $d^*(S)=d^* = \max \frac{|d_k|}{|d_m|}$, where minimum and maximum are taken over all possible combinations of coordinate signs of the solution ${\bf x}(t)$, for {\it any} $0 \le s \le t$.

\medskip
Let there exist positive numbers $M, \beta$ such that
\begin{equation}
e^{-\int_s^t \alpha^*(\tau)\, d\tau} \le M e ^{-\beta(t-s)}, \quad 0 \le s \le t .
\label{upper015}
\end{equation}

\smallskip 

Fixe now an arbitrary $t^* >0$ and consider interval $[0,t^*]$. if our original coefficients are locally integrable, they can be approximated arbitrarily accurately by a continuous functions. In turn, a continuous function can be approximated arbitrarily accurately by an analitic function. As a result, instead of the integrable $B^*(t)$, we obtain an analitic $\bar{B}^*(t)$, such that 
\begin{eqnarray}
\int_0^{t^*} \|{B}^*(\tau)-\bar{B}^*(\tau)\| d\tau  \le \varepsilon.
\label{upper0350}
\end{eqnarray}

Denote now by $W(t,s)$ and $\bar{W}(t,s)$ the Cauchy operators for  (\ref{upper31}) and the respective system with matrix $\bar{B}^*(t)$. Then, if (\ref{upper015}) holds, in accordance with Lemma 3.2.3 \cite{DK} we obtain
\begin{eqnarray}
\|W(t,s)-\bar{W}(t,s)\| \le Md^* e^{-\beta (t-s)}\left(e^{Md^*\int_s^t\|{B}^*(\tau)- \bar{B}^*(\tau)\| d\tau} - 1\right) \nonumber \\
\le  Md^* e^{-\beta (t-s)}\left(e^{Md^*\varepsilon } -1\right).
\label{upper035}
\end{eqnarray}

\smallskip
Hence we have the following statement.

\smallskip

\begin{lemma} Let all $b^*_{ij}(t)$ be locally integrable on $[0,\infty)$. Let inequality (\ref{upper015})  holds. Then 
\begin{eqnarray}
\|{\bf x}\left(t\right)\|_1 \le d^*(S)Me^{-\beta (t-s)}\|{\bf x}(s)\|_{1}, \label{upper32}
\end{eqnarray}
\noindent for any solution of (\ref{upper31}) and any $0 \le s \le t$.
\end{lemma}

\medskip

Let us now return to countable system (\ref{hom11}) and consider the corresponding  truncated system
\begin{equation}
\frac{d}{dt} {\vu}(n,t) = B^*(n,t) {\vu}(n,t), \label{upper11}
\end{equation}
where $B^*(n,t)=\left(b^*_{ij}(t)\right)_{i,j=1}^n$.

Below we will identify the finite vector with entries
$(a_1,\dots, a_n)$ and the infinite vector with the same first
$n$ coordinates and the others equal to zero.

\medskip

Rewrite system (\ref{upper11}) as
\begin{equation}
\frac{d}{dt} {\vu}(n,t) = B^*(t) {\vu}(n,t) + \left( B^*(n,t) -B^*(t) \right) {\vu}(n,t). \label{upper12}
\end{equation}

Denote by $V(t,s)$ and $V(n,t,s)$ the Cauchy operators for (\ref{hom11}) and (\ref{upper11}) respectively.

Suppose that $n > S$, and that, in addition 
\begin{equation}
{\vu}(0)={\vu}(n,0)={\vu}(S,0), \quad \|{\vu}(0)\|_1 \le 1. 
\label{upper125}
\end{equation}

Then one has from (\ref{hom11})
\begin{equation}
{\vu}(t) = V(t){\vu}(0)=V(t){\vu}(n,0).\label{upper13}
\end{equation}

On the other hand, from (\ref{upper12}) we have

\begin{equation}
{\vu}(n,t) = V(t){\vu}(n,0)+ \int_0^t V(t,\tau)\left( B^*(n,\tau) -B^*(\tau) \right) {\vu}(n,\tau)\, d\tau .\label{upper14}
\end{equation}

Hence in any norm we obtain the bound

\begin{equation}
\|{\vu}(t)- {\vu}(n,t)\| \le \int_0^t \|V(t,\tau)\|\|\left( B^*(n,\tau) -B^*(\tau) \right) {\vu}(n,\tau)\|\, d\tau .\label{upper15}
\end{equation}

\smallskip

Denote  $\sup \frac{|d_k|}{|d_m|} = \hat{d} < \infty$, where supremum is taken over all possible combinations of coordinate signs of the solution $\vu(t)$ of (\ref{hom11}), under assumption $|k-m|=1$.

\smallskip
Put now $D^*=diag\left(d^*(1),d^{*}(2),\dots\right)$. 
\smallskip

Note that according to (\ref{b**02}) the matrix $B^{**}(t) $ has nonzero entries only on the main diagonal and at most $R$ diagonals above and below it. Then 
\begin{equation} \|B^{*}(t)\|_{1D^*} =  \|B^{**}(t)\|_{1} \le K=2L\hat{d}^R,\label{upper16}
\end{equation}
\noindent  for almost all $t \ge 0$. 
Then
\begin{equation}
\|V(t,s)\|_{1D^*} \le e^{K(t-s)} \le e^{Kt^*}.
\label{upper17}
\end{equation}
On the other hand, the first $n-R$ columns of the matrix $\left( B^*(n,\tau) -B^*(\tau) \right)$ equil to zero for any $\tau \ge 0$, hence the first $n-R$ coordinates of the corresponding vector $\left( B^*(n,\tau) -B^*(\tau) \right) {\vu}(n,\tau)$ equil to zero and
\begin{equation}
\|\left( B^*(n,\tau) -B^*(\tau) \right) {\vu}(n,\tau)\|_{1D^*} \le K\sum_{k=n-R}^n \hat{d}^k|u_k(t)|.
\label{upper18}
\end{equation}

Put $D^{**}=diag(\hat{d}^{2},\hat{d}^{4},\dots)$ and ${\bf w}^*(t)=D^{**}{\bf u}(t)$.

Then instead of (\ref{upper16}) and (\ref{upper17}) we have
\begin{equation}\|B^*(t)\|_{1D^{**}} = \|D^{**}B^{*}D^{{**}^{-1}}(t)\|_{1} \le K^*=2L\hat{d}^{2R},\label{upper19}
\end{equation}
\noindent  and
\begin{equation}
\|V(t,s)\|_{1D^{**}} \le e^{K^*(t-s)} \le e^{K^*t^*}.
\label{upper20}
\end{equation}
\noindent respectively.

\smallskip

Then
\begin{eqnarray}
\|{\bf u}\left(n,t\right)\|_{1D^{**}} =  \sum_{k = 1}^n \hat{d}^{2k}|u_k(n,t)| \le  e^{K^*t^*}\sum_{k = 1}^S \hat{d}^{2k}|u_k(n,0)|. \label{upper21}
\end{eqnarray}

Then (\ref{upper21}) and (\ref{upper125}) imply the bound

\begin{eqnarray}
\hat{d}^{2n-2R}\sum_{k = n-R}^n |u_k(n,t)| \le  \sum_{k = 1}^n \hat{d}^{2k}|u_k(n,t)| \le  e^{K^*t^*} \hat{d}^{2S}. \label{upper22}
\end{eqnarray}

Then
\begin{eqnarray}
\sum_{k=n-R}^n \hat{d}^k|u_k(n,t)| \le \hat{d}^{n}\sum_{k = n-R}^n |u_k(n,t)| \le  e^{K^*t^*} \hat{d}^{2S+2R-n}.
\label{upper23}
\end{eqnarray}

Finally, for the right-hand side of (\ref{upper15}) we have the bound

\begin{equation}
\int_0^t \|V(t,\tau)\|_{1D^*} \|\left( B^*(n,\tau) -B^*(\tau) \right) {\vu}(n,\tau)\|_{1D^*}\, d\tau \le 
e^{Kt^*}Kt^*e^{K^*t^*} \hat{d}^{2S+2R-n},\label{upper235}
\end{equation}
\noindent which tends to zero at $n \to \infty$.

Hence we have the following statement.

\bigskip

\begin{lemma} Let assumptions of Lemma 1 be fulfilled for any $S$. Then, under assumption (\ref{upper125}), and  for any fixed  $\varepsilon > 0$, $t^*>0$, we get   
$\|{\bf u}(t)-{\bf u}(n,t)\|_{1D^*} < \varepsilon$ for sufficiently large $n$, for any $t \in [0,t^*]$.
\end{lemma}

\bigskip

As a result, Lemmas 1 and 2 guarantee an estimate of the form 
\begin{eqnarray}
\|{\bf u}\left(t\right)\|_1 \le Me^{-\beta t}\|{\bf u}(0)\|_{1D^*}. \label{upper237}
\end{eqnarray}

\smallskip

Consider now two arbitrary solutions $\vp^*(t)$ and $\vp^{**}(t)$ of the forward Kolmogorov system (\ref{ur01}) with the corresponding initial conditions $\vp^*(0)$ and $\vp^{**}(0)$. Denote by $\vp_0^*(t)$ and $\vp_0^{**}(t)$ the respective vector functions with coordinates $1, 2, \dots$ (i.e. without zero coordinates).

One can write $\vu(t)=T\left(\vp_0^*(t)-\vp_0^{**}(t)\right)$. Then (see for instance \cite{Zeifman2014c}) the following inequality holds: $\|\vp^*(t)-\vp^{**}(t)\|_1 \le \frac{2}{d} \|\vu(t)\|_1$.

\bigskip

Finally we obtain the following statement.

\begin{theorem}{}
Let the assumptions of Lemma 1 hold for any natural $S$. Then $X(t)$ is weakly ergodic and the following bound on the rate of convergence holds:
\begin{eqnarray}
\|\vp^*(t)-\vp^{**}(t)\|_1\le \frac{2M}{d}e^{-\beta t}\|\vp_0^*(0)-\vp_0^{**}(0)\|_{1D^*}. \label{upper238}
\end{eqnarray}

\end{theorem}

\medskip

\begin{remark} A  specific model (which belongs to both  classes $\rm (iii)$ and $\rm (iv)$) was investigated in \cite{Satin2019} by the method described here. 
\end{remark}

\medskip

Namely, in this paper, the queueing model with possible transitions and respective intensities  of single arrival $\la(t)$ and service of group of two customers $\mu(t)$  was considered. Hence
$$B^*(t)=$$
\begin{eqnarray}
 = { \left(
\begin{array}{ccccccc}
-\lambda(t)   & -\mu(t)                      & \mu(t)                       & 0  & 0 & 0 &\cdots\\
 \lambda(t)   & -\left(\lambda(t)+\mu(t)\right) & 0                         & \mu(t)  & 0 & 0 &\cdots \\
  0        & \lambda(t)                   & -\left(\lambda(t)+\mu(t)\right) & 0 & \mu(t)  & 0 & \cdots \\
  0        & 0                         & \lambda(t)                   & -\left(\lambda(t)+\mu(t)\right)  & 0 & \mu(t) & \cdots\\
  0        & 0                         & 0                         & \lambda(t) & -\left(\lambda(t)+\mu(t) \right)  & 0 & \cdots\\
\ddots & \ddots & \ddots & \ddots & \ddots & \ddots  & \ddots \\
\cdots & \cdots & \cdots & \cdots & \cdots & \cdots & \cdots
\end{array}
\right)} .\nonumber
\end{eqnarray}

\smallskip

Let $\delta >1$ be a positive number. Put 

$d_1=1$, $d_2=1/\delta$, $d_k=\delta^{k-2}$, $k \ge 3$ if all coordinates of solutions are positive;

 $|d_k|=\delta^{k-1}$, $k \ge 1$ otherwise.

Then one has, 
\begin{eqnarray}
\alpha^*(t) \ge \min \left[\la(t)\left(1-\delta^{-1}\right), \mu(t)
\left(1+\delta\right)-\right. \nonumber
\\ \left. \la(t)\left(\delta^{2}-1\right), \mu(t)
\left(1-\delta^{-1}\right)- \la(t)\left(\delta-1\right) \right].
\label{finbound01}
\end{eqnarray}
\noindent   Moreover $d = \delta^{-1}$, $\hat{d}=\delta$, $d^*_k = \delta^{k-1}$, for $k \ge 1$.

In particular, if the process $X(t)$ is homogeneous i.e.,  $\la(t)=\la$ and $\mu(t)=\mu$ are positive numbers, then $\int_0^{\infty}\alpha^*(t) dt = +\infty$  is equivalent to $\alpha^* >0$ and this
is equivalent to $0 < \la < \mu$. Put $\delta= \sqrt{\frac{\mu}{\la}}$. Hence,
\begin{equation}\alpha^* = \min\left[\left(\sqrt{\mu} -\sqrt{\la}\right)^2,
\la\left(1-\sqrt{\frac{\la}{\mu}}\right)\right].\label{finbound05}
\end{equation}

\clearpage

\section{Acknowledgment}
This research was supported by Russian Science Foundation under grant 19-11-00020.

\end{document}